\documentclass[12pt]{article}

\usepackage{amsmath,amssymb,latexsym,float,epsfig,subfigure}

\usepackage{latexsym}
\usepackage{amsthm}
\usepackage{cases}
\usepackage{epsfig}
\usepackage{subfigure}
\usepackage{scalefnt}
\usepackage{booktabs}
\usepackage[english]{babel}
\usepackage{graphicx}
\usepackage{titlesec}
\usepackage{epsfig}
\usepackage{enumerate}
\usepackage{caption}
\usepackage{subfigure}
\usepackage{booktabs}
\usepackage[english]{babel} 
\usepackage{graphicx}
\usepackage{multirow}
\usepackage{rotating}
\usepackage{comment}
\usepackage{algorithm}
\usepackage{algpseudocode}

\usepackage{algpseudocode}
\usepackage{multirow}

\usepackage{xargs}                      
\usepackage[pdftex,dvipsnames]{xcolor}  
\usepackage[colorinlistoftodos,prependcaption,textsize=tiny]{todonotes}
\newcommandx{\unsure}[2][1=]{\todo[linecolor=red,backgroundcolor=red!25,bordercolor=red,#1]{#2}}
\newcommandx{\change}[2][1=]{\todo[linecolor=blue,backgroundcolor=blue!25,bordercolor=blue,#1]{#2}}
\newcommandx{\info}[2][1=]{\todo[linecolor=OliveGreen,backgroundcolor=OliveGreen!25,bordercolor=OliveGreen,#1]{#2}}
\newcommandx{\improvement}[2][1=]{\todo[linecolor=Plum,backgroundcolor=Plum!25,bordercolor=Plum,#1]{#2}}
\newcommandx{\thiswillnotshow}[2][1=]{\todo[disable,#1]{#2}}

\usepackage{fancyhdr}
\usepackage{setspace}
\usepackage[section]{placeins}
\usepackage{xcolor}
\usepackage{amsmath}
\makeatletter
\renewcommand{\Function}[2]{%
  \csname ALG@cmd@\ALG@L @Function\endcsname{#1}{#2}%
  \def\jayden@currentfunction{#1}%
}
\newcommand{\funclabel}[1]{%
  \@bsphack
  \protected@write\@auxout{}{%
    \string\newlabel{#1}{{\jayden@currentfunction}{\thepage}}%
  }%
  \@esphack
}
\makeatother

\usepackage{breqn}
\usepackage{gensymb}
\usepackage{array, makecell}
\usepackage{pseudocode}

\usepackage{url} 
\usepackage{appendix}
\usepackage{ragged2e}
\usepackage{physics}
\usepackage{natbib} 
\usepackage[thinlines]{easytable}

\newtheorem{proposition}{Proposition}

\newtheorem{lemma}{Lemma}

\usepackage[nocomma]{optidef}
\usepackage{mathtools,zref-savepos}
\newtagform{roman}[]()

\topmargin by -\headsep \textheight 8.9in \oddsidemargin 0pt
\evensidemargin \oddsidemargin \marginparwidth 0.5in \textwidth
6.5in
\date{}
\usepackage{graphicx}
\usepackage{bbm}
\usepackage{mathtools}
\usepackage{graphicx}
\usepackage{caption}
\usepackage{subcaption}

\usepackage{algpseudocode}
\allowdisplaybreaks
\usepackage{mathrsfs}
\usepackage{tikz}
\usetikzlibrary{decorations.pathreplacing}
\usetikzlibrary{backgrounds}
\usetikzlibrary{arrows,shapes,positioning}

\usepackage{calligra}

\usepackage{authblk}

\providecommand{\keywords}[1]
{
  \small
  \textbf{\textit{Keywords---}} #1
}

\title{Approximate Solutions for Multi-Trip Route Planning in Time-Sensitive Situations}
\author[1]{Bahar \c{C}avdar\thanks{cavdab2@rpi.edu}}
\author[2]{Joseph Geunes\thanks{geunes@tamu.edu}}
\author[3,2]{Xiaofeng Nie\thanks{xiaofengnie@tamu.edu}}
\author[4]{Yue Wang\thanks{ywang23@arizona.edu}}

\affil[1]{\small{Department of Industrial and Systems Engineering, Rensselaer Polytechnic Institute, Troy, NY }}
\affil[2]{Wm Michael Barnes '64 Department of Industrial and Systems Engineering, Texas A\&M University, College Station, TX}%
\affil[3]{Department of Engineering Technology and Industrial Distribution, Texas A\&M University, College Station, TX}
\affil[4]{Department of Systems and Industrial Engineering, The University of Arizona, Tucson, AZ}

\begin{document}

\maketitle

\begin{abstract}
We consider emergent situations that require transporting individuals from their locations to a facility using a single capacitated vehicle, where transportation duration has a negative impact on the individuals. A dispatcher determines routes to maximize total satisfaction. We call this problem the Ambulance Bus Routing Problem.  We develop efficient approximate policies for the dispatcher to allocate individuals to multiple routes, characterize an optimal solution of the relaxed approximate model, and devise a heuristic to obtain a near-optimal integer solution quickly.
\end{abstract}


\section{Introduction}
\label{Sec:Intro}
The majority of routing problems consider how to determine vehicle routes to visit demand locations either to pick up an item or to make a delivery. Most of these problems are solved using efficiency-based objectives such as minimizing the cost of operations, measured in terms of distance, time, or fuel cost, without considering the impact on the individuals served by the routes. The classical Traveling Salesman Problem (TSP) and Vehicle Routing Problem (VRP) fall within this group.  Other problems, such as the Traveling Repairman Problem (TRP), use a more customer-centric objective function and consider cumulative waiting times of all customers, reflecting the negative impact of delays in service delivery from customers' perspectives. In these problems, it is also assumed that the service for a customer is considered to be completed upon arrival at the customer location. However, many routing applications exist where service requires transporting customers from their origin to a specific destination, any delays in transportation have a negative impact on the customers, and the decisions need to be made in real time.

We study a logistics problem where geographically dispersed individuals need transportation to a certain location, with their condition deteriorating over time. Such problems appear in many areas, including humanitarian logistics, military operations, and service systems. For example, in preparation for a natural disaster, e.g., a hurricane, residents in the affected areas may need to be transported to safe locations. Over time, the risk of being negatively impacted by the hurricane increases. In addition, in the aftermath of a disaster, e.g., an earthquake, there may be a large number of injured individuals spread across a region with different levels of severity. Emergency vehicles, such as ambulance buses, can be used to transport these individuals to a healthcare facility. Group routing decisions in humanitarian logistics become especially critical in urban areas because of the large number of transit-dependent individuals. For example, in New York City, 47\% of the people do not own or have access to a car \citep{NYC}. In military operations, in  a developing or escalating situation, personnel in a region may need to be transferred to a safe location, with  any delay increasing the risk to the personnel. Finally, we observe similar phenomena in, for example, airport parking services. In these systems, after customers park in a lot, they wait for shuttle transportation to the terminal. The shuttle operator determines which customers will be transported on each shuttle run. While picking up more customers in each run will benefit those boarding last, it can result in delays for passengers already on the shuttle, leading to increased dissatisfaction due to the extended travel time to the terminal.

Motivated by these applications, we study efficient routing policies for group transportation under emergent or urgent situations. We call the underlying routing problem the Ambulance Bus Routing Problem (ABRP).
We form a connection between the ABRP and the TSP, although the ABRP is computationally more demanding due to the nature of the cumulative objective function and additional route-assignment decisions.
This connection allows us to develop an approximate model for the routing problem under emergent situations, with two significant benefits. The approximation drives insights for dispatchers to support real-time decision making and also has an interesting structure as a standalone model. 
In our approximate model, we use the asymptotic TSP tour length estimation model from \citet{beardwood1959shortest} and provide structural analysis of optimal solutions when customers are uniformly distributed throughout a region. The resulting problem is, to the best of our knowledge, a new class of equality-constrained knapsack problems with a nonlinear objective containing an interesting special structure. Despite the problem's lack of convexity, our analysis reveals a surprisingly simple and counter-intuitive policy for solving the continuous relaxation of the approximate problem.
We then devise a heuristic for obtaining an integer solution based on these results. In our computational experiments, our heuristic provides high-quality solutions comparable to those obtained by the global optimization solver BARON but in a significantly shorter amount of time.

The rest of the paper is organized as follows. In Section \ref{Sec:Lit}, we review the related literature and position our paper with respect to related routing problems. Section \ref{Sec:Def} presents a problem definition along with a mathematical programming formulation for the ABRP to determine exact optimal routes for completeness.
In Section \ref{Sec:Soln}, we present our approximate model and analyze the characteristics of optimal routing strategies for both uncapacitated and capacitated versions under a linear satisfaction function. We present and discuss our computational experiment results in Section \ref{Sec:Comp}. In Section \ref{Sec:Conc}, we provide a summary of our work and present some future research directions.

\section{Literature Review}
\label{Sec:Lit}
This paper considers routing decisions in emergent or urgent situations, with prevalent applications in disaster preparation and response, and humanitarian logistics. Therefore, we review two main streams of research. We first review classical routing problems, such as the TSP, the VRP, and their variants. Due to the extensive nature of this research stream, we do not attempt to provide a comprehensive review, but rather connect our underlying routing problem to the literature. In the second part, we review studies on routing problems specifically in emergency response. 

\emph{Classical Routing Problems:} We classify routing problems where the objective consists of minimizing travel cost within this group. The TSP and the VRP are the two most popular examples in this group. The body of research on classical routing problems is vast and well-established. \citet{laporte2009fifty} is one of the most notable reviews on the VRP focusing on different solution methods. \citet{eksioglu2009vehicle} and \citet{braekers2016vehicle} provide a taxonomic review for the VRP. \citet{cook2011traveling} focus on large-scale TSPs and review solution methods.

The TSP and the VRP also have numerous variants. \citet{Vidal2020-tw} review new variants of the VRP focusing on new performance metrics and objectives, how to integrate route planning with other business decisions, and the need for more precise modeling approaches. Some variants of classical routing problems have emerged to accommodate the requirements of business environments. For example, in multi-trip TSPs and VRPs, vehicles can operate on multiple consecutive routes \citep{cattaruzza2016vehicle}. On the other hand, other variants use different objectives that are more customer-centric. Cumulative TSPs and VRPs, where the objective is to minimize the total service time for each customer, fall within this group. \citet{corona2022vehicle} review routing problems with cumulative objectives emphasizing  applications in healthcare, disaster relief operations, and maintenance. Despite the high-impact application areas, the literature on cumulative routing problems is still sparse. Our problem involves a multi-trip routing decision with a cumulative objective function where we calculate the cumulative negative impact of service duration across all customers.

\emph{Emergency Response Routing:} Routing decisions also arise in evacuation route planning. In this stream, most of the studies focus on car-based evacuation planning, where individuals are expected to follow planned routes either in preparation for or in the aftermath of a disaster \citep{Esposito_Amideo2019-jj}. The importance of effective group routing planning was reinforced after Hurricane Katrina (2005), during which long evacuation delays were observed due to traffic congestion \citep{Lakshay2019-yf}. Hence, there has been a growing effort to develop evacuation route strategies using large-capacity vehicles. \citet{Bish2011-sr} provides one of the earliest works focusing on group routing for transit-dependent populations in preparation for an approaching threat. Given a set of nodes, including the current locations of the buses, the evacuees, and the shelters, the objective is to determine routes with the minimum total evacuation time. This problem is called the Bus Evacuation Problem (BEP). The paper by \citet{Goerigk2013-au} builds upon the BEP and focuses on developing efficient solution methodologies. The authors propose new lower and upper bounds for the problem and incorporate these in a branch-and-bound framework. \citet{Goerigk2014-mz} study the BEP under the assumption that the exact number of evacuees is not known in advance. Their primary goal is to decide when to dispatch a bus to minimize the maximum travel time over all buses while transporting all evacuees to shelters. The nature of this type of objective function implicitly enforces fairness although it is still efficiency-based. \citet{Swamy2017-of}, \citet{Lakshay2019-yf}, and \citet{Zhao2020-tf} provide additional studies of bus routing problems in emergency response with a goal of minimizing the total travel time.

\citet{Dikas2016-fq} emphasize that solving bus routing problems that minimize total travel time is not the most appropriate approach for applications in humanitarian logistics and healthcare. Their paper builds on the work by \citet{Bish2011-sr} and studies routing decisions that minimize the total time required to transport casualties to medical facilities, which is equivalent to minimizing the average transportation time per individual. They call this version of the problem the Casualty Evacuation Problem. Their computational experiments reveal the inherent challenges in solving this kind of routing problem; they allow up to 6 hours of computation time for instances with 8 nodes, and 24 hours for instances with 15 nodes. \citet{Jin2015-dg} study the problem of transporting patients to healthcare centers in the aftermath of a disaster. They consider patient survival probabilities and seek routes that maximize the expected number of survivors. 

Our work is closely related to these studies on emergency response routing, but incorporates two distinguishing features. In contrast to the majority of prior works, our paper considers the negative impact of service duration on individuals. In particular, our objective function accounts for the decreasing satisfaction or wellness of the individuals as transportation duration increases. In addition, motivated by the time-sensitive nature of emergency routing decisions, we present an approximation model to compute general policies for constructing patient transportation routes rather than determining exact route sequences.

\section{Problem Description}
\label{Sec:Def}
In this section, we first describe the underlying problem environment and definition. 
We use the disaster response setting to describe the problem, although our solution methods can be adapted to other applications in service systems and military operations.
Since our approximate model incorporates the TSP tour length estimation model from \cite{beardwood1959shortest} into the ABRP, which determines exact routes, we first present the ABRP formulation for completeness.

We consider an emergent situation in the aftermath of a disaster where a number of affected people are dispersed throughout a geographical area. A rescue team manages a single ambulance bus with a capacity of carrying $C$ individuals that picks up affected individuals and transports them to a central facility. We denote the underlying undirected network by $G=(V, E)$, where $V$ is the set of nodes and $E$ is the set of edges. Let $N$ denote the number of affected individuals dispersed throughout the area who require transport to the central facility. For convenience we use the terms individuals and nodes interchangeably. Node $i$ is contained within $V=\{0, 1, \dots, N\}$, where node 0 represents the central facility and $d_{ij}$ is the travel time between nodes $i$ and $j \neq i$ for $(i,j)\in E$.

We assume that an individual's service is completed upon arrival at the central facility and that an affected individual's condition deteriorates with time. In an emergency response or humanitarian logistics setting, an individual's condition may correspond to a survival probability, while in other contexts it may correspond to some measure of service quality. 
To quantify this, we use a function $S(t)$ that calculates the status of an individual at time $t$, where $S(t)$ is a non-increasing function of $t$. In service applications, $S(t)$ may measure customer satisfaction (e.g., in airport parking). We will generically refer to $S(t)$ as the satisfaction level.  For ease of exposition and analytical tractability we assume that $S(t)$ is a linearly decreasing function of $t$, i.e., $S(t)=a-bt$. For some empirical survival probability functions in emergency response applications, please refer to \citet{erkut2008ambulance} and \citet{mills2013resource}.

In many of the applications we have discussed, in the interest of fairness, the service sequence should not violate the sequence in which service calls are received. To model this, without loss of generality, we number customer node indices based on  the order in which service requests are received. That is, $j>i$ implies that the service request from node $i$ arrived before that of node $j$. Accordingly, if individual $i$ is transported on route $k$ and individual $j$ is transported on route $k'$, then $k\leq k'$ for any $i<j$.

The problem requires determining the sequential routes for the service vehicle that maximize the total satisfaction level of all individuals requiring transportation at the time of service completion. We call this problem the ABRP. We present a mixed-integer programming formulation for the ABRP below.

\textbf{Decision variables}
\begin{itemize}
    \item $y_{ik}$: 1 if node $i$ is on route $k \in \{1, \dots, N\}$, 0 otherwise.
    \item $x_{ij}^{k}$: 1 if the ambulance bus travels directly from node $i$ to node $j \neq i$ on route $k$, 0 otherwise.
    \item $A_i$: individual $i$'s arrival time at the destination (service completion time).
    \item $t_k$: the total travel time of the $k^{\text{th}}$ route.
\end{itemize}
\begin{align}
\textbf{[P1]} \qquad&\text{Max} && \sum_{i=1}^{N}S(A_i) \label{Mod:obj}  \\
& \text{s.t.} && \sum_{k=1}^{N}y_{ik}=1, && \forall i  \neq 0,\label{Mod:Mustvisit}\\
&&& \sum_{i=1}^{N} y_{ik} \leq C, && \forall k, \label{Mod:Cap}\\
&&& t_k = \sum_{(i,j)\in E} d_{ij} x_{ij}^{k}, && \forall k, \label{Mod:Tour}\\
&&& A_i \geq \sum_{h=1}^{k} t_{h} - M (1- y_{ik}), && \forall i \neq 0, k, \label{Mod:Arrivaltime}\\
&&& y_{ik}= \sum_{j \in V \setminus \{i\}} x_{ij}^{k} = \sum_{j \in V \setminus \{i\}}x_{ji}^{k}, && \forall i, ,k  \label{Mod:flow}\\
&&& u_0 = 1, \label{Mod: mtz_subtour_1}\\
&&& 2 \leq u_i \leq N+1 ,  && \forall i \neq 0,  \label{Mod: mtz_subtour_2}\\
&&& u_j - u_i \geq 1-N(1-x_{ij}^k), && \forall i\neq 0, j \neq 0, j \neq i, k, \label{Mod: mtz_subtour_3}\\
&&& y_{ik} \in \{0,1\}, && \forall i, k, \\
&&& x_{ij}^{k} \in \{0,1\}, && \forall i, j \neq i, k.
\end{align}
\normalsize

The objective function (\ref{Mod:obj}) maximizes the cumulative satisfaction of all individuals at each individual's service completion time. Constraint set (\ref{Mod:Mustvisit}) ensures that each service request is met, while constraint set (\ref{Mod:Cap}) ensures that the ambulance bus capacity is not violated. Constraints (\ref{Mod:Tour}) calculate the total travel time on each route $k$. Constraint set (\ref{Mod:Arrivaltime}) calculates the service complete time of individual $i$ using a big $M$. Constraint set (\ref{Mod:flow}) ensures the connectivity of each route and constraints (\ref{Mod: mtz_subtour_1})-(\ref{Mod: mtz_subtour_3}) eliminate subtours within each route \citep{miller_integer_1960}.

\begin{proposition}
\label{Prop:NP-complete}
The ABRP is NP-complete. 
\end{proposition}
The proof is in Appendix \ref{Ap:NPComplete}.

Note that a significant difference between the ABRP and the TRP results from the time at which service is completed. In the ABRP, we assume that service is completed when the individual arrives at the final destination. If $T_k$ denotes the time when route $k$ is completed, i.e., $T_k=\sum_{h=1}^{k}t_h$, then the satisfaction level for all customers on route $k$ equals $S(T_k)$. 
On the other hand, the corresponding value is customer-specific in the TRP.

Although \textbf{[P1]} is an NP-complete problem, the cumulative nature of the objective function in the ABRP permits characterizing the following key structural property of optimal solutions.

\begin{proposition}\label{prop:routeorder}
Consider a solution to an ABRP containing $K$ routes, with $n_k$ individuals on route $k$. 
If $n_{k} \leq n_{h}$ and $t_k > t_h$, then any solution in which route $k$ is completed before route $h$ is suboptimal.
\end{proposition}

The proof is in Appendix \ref{Ap:routeorder}. Proposition \ref{prop:routeorder} implies that a candidate route with a larger number of individuals and faster execution time should be completed prior to routes with longer duration and fewer individuals. Making use of this property, however, requires comparing pairs of candidate routes, of which there is an exponential number in $N$.  Because of this, finding good solutions for the ABRP in real time, especially under time-sensitive conditions, is extraordinarily challenging. 
On the other hand, the characteristics of the ABRP create opportunities for an approximation model, which leads to both very fast solutions and insights, supporting route-planners' real-time decision-making in time-sensitive environments. 
We will therefore focus on some key structural properties (similar to Proposition \ref{prop:routeorder}) of optimal solutions for an approximation of \textbf{[P1]} to devise efficient algorithms.
We present our approximate model and solution approach in the next section.

\section{Approximation Using TSP Tour Length Estimation}
\label{Sec:Soln}

In an optimal solution to the ABRP, for any given allocation of nodes to routes, the pickup sequence within each route must follow an optimal TSP route. Using this observation, we replace the cost of each route and the routing decisions in the ABRP formulation with an estimation of the optimal TSP tour length to simplify the allocation decisions.
After determining an allocation of nodes to routes, we can solve a TSP within each route. 

There are several tour length estimation models targeting different attributes of an underlying network \citep{cavdar2015distribution}. Among these, we use the model from \citet{beardwood1959shortest} due to its simplicity and high estimation power. In their seminal work, \citet{beardwood1959shortest} showed that the optimal TSP tour length or duration $T$ for serving $n$ customers distributed over an area $A$ can be approximated by $T\sim \beta \sqrt{nA}$,
where $\beta$ is a constant that depends on the dispersion of the nodes, computationally verified to be around 0.72 \citep{cook2011traveling} for uniform node dispersion. While this model is asymptotic and performs better as $n$ increases, recent studies have shown that an appropriate $\beta$ value can be computed for a small number of nodes as well \citep{vinel2018probability}.

Recall that in the ABRP, the order of the service should follow the order of the requests. Therefore, we assume that any subset of $n_k$ nodes is also uniformly distributed on the region of area $A$. By letting $\kappa=\beta \sqrt{A}$, the time at the completion of the $k^{\text{th}}$ route equals
$T_k=\kappa \sum_{i=1}^{k} \sqrt{n_i}$.
For a given $k \in\{1,2, \dots, N\}$, our goal then is to solve
\begin{align}
\textbf{[P2]} \qquad&\text{Max} && \sum_{i=1}^k n_i S\left(\kappa \sum_{j=1}^i \sqrt{n_j}\right) \nonumber\\
& \text{s.t.} &&\sum_{i=1}^{k} n_i = N,\nonumber\\
&&&n_i \leq C,  && i=1, \dots, k,\nonumber \\
&&& n_i \in \mathbb{Z}^+, && i = 1, \ldots, k. \nonumber
\end{align}
\normalsize

The objective function for a linear satisfaction function $S(t)=a-bt$ can be written as follows:
\begin{equation}
\text{Max} \qquad aN - b \kappa \sum_{i=1}^{k}n_i \bigg( \sum_{j=1}^{i} \sqrt{n_j}\bigg), \label{eq:P3_maxObj}
\end{equation}
where $\sum_{i=1}^{k}n_i=N$ and each $n_i$ is a nonnegative integer. The corresponding optimal solution for \textbf{[P2]} can thus be obtained by solving
\begin{eqnarray}
\textbf{[P3]} \qquad \qquad\mbox{Min} & \sum_{i=1}^k n_i \left(\sum_{j=1}^i \sqrt{n_j}\right)   & \nonumber  \\
\mbox{s.t.} & \sum_{i=1}^k n_i = N, & \nonumber \\
& n_i \leq C, & i = 1, \ldots, k, \nonumber\\
& n_i \in \mathbb{Z}^+,  & i = 1, \ldots, k. \nonumber
\end{eqnarray}
\normalsize

We call this problem the Approximate Ambulance Bus Routing Problem (AABRP), which corresponds to an equality-constrained version of the continuous knapsack problem with a nonlinear objective. In the general version of the AABRP, we have $k=N$ routes.

\begin{proposition}
\label{Prop:ordering}
In an optimal solution to the AABRP in the form of \textbf{[P3]} where $n_i>0$, $n_{i+l}>0$, and $l$ is a positive integer, we have $n_i \geq n_{i+l}$.
\end{proposition}

The proof of the proposition is in Appendix \ref{Ap:ordering}.

The following subsection presents our solution method development for the AABRP. We first consider the uncapacitated problem, and then extend our results to the capacitated problem.

\subsection{Uncapacitated AABRP} 
Our solution method for the AABRP is based on the characteristics of the optimal solution of its relaxation. We first consider the continuous relaxation of the uncapacitated problem:
\begin{eqnarray}
\textbf{[P4]} \qquad\mbox{Min} & \sum_{i=1}^k n_i \left(\sum_{j=1}^i \sqrt{n_j}\right)   &  \nonumber\\
\mbox{s.t.} & \sum_{i=1}^k n_i = N, & \nonumber \\
& -n_i \le 0,  & i = 1, \ldots, k. \nonumber
\end{eqnarray}
\normalsize

We claim that all decision variable values will be strictly positive in an optimal solution. 
\begin{proposition}
\label{Prop:Strictpositve}
In an optimal solution of \textbf{[P4]}, $n_i > 0$ for all $i$. 
\end{proposition}

The proof is in Appendix \ref{Ap:Strictpositve}.

We can show that the above objective function is neither convex nor concave, and because the constraints form a linear feasible region, the Karush–Kuhn–Tucker (KKT) conditions are necessary, but not sufficient for optimality. Despite this, our analysis will permit identifying a unique KKT point that satisfies the conditions of Propositions \ref{Prop:ordering} and \ref{Prop:Strictpositve}. Let $\alpha$ be the KKT multiplier associated with the equality constraint, and let $\lambda_i$ denote a KKT multiplier associated with the $i^{\text{th}}$ nonnegativity constraint. The KKT conditions may be written as 
\begin{eqnarray}
 & \frac{3\sqrt{n_i}}{2} + \sum_{j=1}^{i-1}\sqrt{n_j} + \frac{\sum_{j=i+1}^{k}n_j}{2\sqrt{n_i}} - \alpha - \lambda_i = 0,  &  i = 1, \ldots, k, \label{eq:KKT1}  \\
 & \sum_{i=1}^k n_i = N, & \\
 & \lambda_i n_i = 0, & i = 1, \ldots, k,\\
 & n_i \ge 0, & i = 1, \ldots, k, \\
 & \lambda_i \ge 0, & i = 1, \ldots, k.
\end{eqnarray}

Observe that $n_i>0$ implies that $\lambda_i=0$; thus Proposition \ref{Prop:Strictpositve} implies that 
\begin{eqnarray}
 & \frac{3\sqrt{n_i}}{2}  + \sum_{j=1}^{i-1}\sqrt{n_j} + \frac{\sum_{j=i+1}^{k}n_j}{2\sqrt{n_i}} = \alpha, &\label{Eq:Must}
\end{eqnarray}
for $i=1,\ldots,k$, while further implying that for any $i=1,\ldots,k$ and $\ell\ne i$ we have 
\begin{eqnarray}
 & \frac{3\sqrt{n_i}}{2}  + \sum_{j=1}^{i-1}\sqrt{n_j} + \frac{\sum_{j=i+1}^{k}n_j}{2\sqrt{n_i}} = \frac{3\sqrt{n_{\ell}}}{2}  + \sum_{j=1}^{\ell-1}\sqrt{n_j}  + \frac{\sum_{j=\ell+1}^{k}n_j}{2\sqrt{n_{\ell}}}   . & \nonumber\label{Eq:altKKT}
\end{eqnarray}

When $i=k$ and $\ell=k-1$ this becomes
\begin{eqnarray}
& \frac{3\sqrt{n_k}}{2}  + \sqrt{n_{k-1}} = \frac{3\sqrt{n_{k-1}}}{2}   + \frac{n_k}{2\sqrt{n_{k-1}}}. & \label{Eq:SolveMust}
\end{eqnarray}
Since any optimal solution with $n_{k-i}>0$ and $n_{k-i-1} >0$ for $i<k-1$ satisfies $n_{k-i-1}\geq n_{k-i}$ by Proposition \ref{Prop:ordering}, a candidate for an optimal point must satisfy 
\begin{eqnarray}
&    \sqrt{n_{k-1}} =  \frac{3+\sqrt{5}}{2}\sqrt{n_{k}}. & \label{Eq:ConsecRel}
\end{eqnarray}

Note that the inverse of $\frac{3+\sqrt{5}}{2}$ is $\frac{3-\sqrt{5}}{2}$ and if we let $\varphi$ denote the golden ratio (GR), i.e., $\varphi=\frac{1+\sqrt{5}}{2}$, we can then write $\frac{3+\sqrt{5}}{2}=1+\varphi$ and $\frac{3-\sqrt{5}}{2}=2-\varphi$. This implies that in an optimal solution to \textbf{[P4]} with $k=2$ routes, the ratio of the number of individuals on the longer tour to that of the shorter tour equals $1+\varphi$, which has an inverse of $2-\varphi$. This relationship also extends to the tour length relationships.

\textbf{Two-route solution:} For illustration purposes, we next consider the case where we have $k=2$ route variables, that is, $n_1+n_2=N$. Following the expression in Equation (\ref{Eq:ConsecRel}), we can write 
$\sqrt{n_{2}} =  \left(2-\varphi \right) \sqrt{n_{1}}$, 
which implies
$$n_1^*  =  \frac{1}{1+\left(2-\varphi \right)^2}N =  \frac{3+\sqrt{5}}{6} N \mbox{ and }
n_2^*  =  \frac{\left(2-\varphi \right)^2}{1+\left(2-\varphi \right)^2} N = \frac{3-\sqrt{5}}{6}N. $$

The corresponding optimal objective function value is calculated as
$\frac{ 1 +  \left(2-\varphi \right)^2 + \left(2-\varphi \right)^3}{\left( 1+\left(2-\varphi \right)^2\right)^{3/2}}N^{3/2}$.

\textbf{Multiple-route solution:} Next, consider cases with $k>2$ route variables, i.e., $2<k\le N$. Then, using Equation (\ref{Eq:Must}) we can write
\begin{eqnarray}
\label{Eq:MultMust}
 & \frac{3\sqrt{n_{k-2}}}{2}   + \frac{n_{k-1}+n_{k}}{2\sqrt{n_{k-2}}} = \frac{3\sqrt{n_{k-1}}}{2}  + \sqrt{n_{k-2}}  + \frac{n_{k}}{2\sqrt{n_{k-1}}}  . &
\end{eqnarray}
We know that an optimal KKT point must satisfy  $\sqrt{n_k}=(2-\varphi)\sqrt{n_{k-1}}$ by Equation (\ref{Eq:ConsecRel}). Therefore, the above becomes
\begin{eqnarray}
 & \frac{\sqrt{n_{k-2}}}{2}   + \frac{(1+(2-\varphi)^2)n_{k-1}}{2\sqrt{n_{k-2}}} =  \frac{(3+(2-\varphi)^2)\sqrt{n_{k-1}}}{2}   , & \nonumber 
\end{eqnarray}
which is equivalent to
\begin{eqnarray}
 & n_{k-2}  - (3+(2-\varphi)^2)\sqrt{n_{k-1}}\sqrt{n_{k-2}}+ (1+(2-\varphi)^2)n_{k-1} = 0. & \label{Eq:MustEqual} 
\end{eqnarray}

Noting that $(2-\varphi)^2 = 5-3\varphi$, we can rewrite Equation (\ref{Eq:MustEqual}) as
\begin{eqnarray}
 & n_{k-2}  - (8-3\varphi)\sqrt{n_{k-1}}\sqrt{n_{k-2}}+ (6-3\varphi)n_{k-1} = 0. & \nonumber 
\end{eqnarray}
This quadratic equation is solved at
\begin{eqnarray}
  \sqrt{n_{k-2}} 
 &  = & \sqrt{n_{k-1}}\left( \frac{8-3\varphi}{2} \pm \frac{\sqrt{49 -27\varphi}}{2} \right). \nonumber 
\end{eqnarray}
Using Proposition \ref{Prop:ordering}, the candidate for an optimal point is
\begin{eqnarray}
  n_{k-2} &   = & n_{k-1}\left( \frac{8-3\varphi}{2} + \frac{\sqrt{49 -27\varphi}}{2} \right)^2. \label{Eq:k-2}
\end{eqnarray}

Let $\nu_{k-i}^k$ denote the ratio between $n_{k-i}$ and $n_{k-i+1}$ when there are $k$ route variables.  
Noting that by Proposition \ref{Prop:Strictpositve}, $n_j>0$ for $j=1,\ldots,k$, and we can write Equation \eqref{Eq:ConsecRel} as  
\begin{equation}
n_{k-1}=\underbrace{(1+\varphi)^2}_{\nu_{k-1}^k}n_k. \label{Eq:NU1}
\end{equation}
Similarly, Equation (\ref{Eq:k-2}), which determines $n_{k-2}$ based on $n_{k-1}$, can be rewritten as 
\begin{equation}
  n_{k-2}  = \underbrace{\left(\frac{3 + \frac{1}{\nu_{k-1}^k} +\sqrt{\left(3+\frac{1}{\nu_{k-1}^k}\right)^2 - 4 \left(1+\frac{1}{\nu_{k-1}^k}\right)} }{2} \right)^2}_{\nu_{k-2}^k}  n_{k-1}. \label{Eq:NU2}
\end{equation}
Continuing in the same fashion recursively, we observe a pattern for the relationship between $n_{k-l}$ and $n_{k-l+1}$ for $l > 2$, i.e., $n_{k-l}= \nu_{k-l}^k n_{k-l+1}$, where
\begin{equation}
\nu_{k-l}^k = \left( \frac{3+\sum_{i=1}^{l-1} \rho_{i}^k + \sqrt{ \left( 3+\sum_{i=1}^{l-1} \rho_{i}^k\right)^2  - 4 \left(   1+\sum_{i=1}^{l-1} \rho_{i}^k \right)   }   }   {2}\right)^2, \label{Eq:NUs}
\end{equation}
where $\rho_{i}^k= \frac{1}{\prod_{j=1}^{i}\nu_{k-j}^k}$ for $i=1,\ldots,k-1$. We next define $\hat{\rho}_i^k=\frac{1}{\rho_{i-1}^k\prod_{j=1}^{k-1}\nu_j^k}=\prod_{j=1}^{k-i}\frac{1}{\nu_j^k}$.


Using Equations \eqref{Eq:NU1} -- \eqref{Eq:NUs} and $\sum_{i=1}^{k} n_i = N$, we obtain the following proposition. 

\begin{proposition}
An optimal solution to \textbf{[P4]} is characterized as 
\begin{equation} 
n_1^*=\frac{N}{1+\sum_{i=1}^{k-1}\hat{\rho}_i^k}  \label{eq:n1}
\end{equation}
and   
\begin{equation}
n_j^*=\hat{\rho}_{k-j+1}^k n_1^* 
\end{equation}
for $j \in \{2,3,\ldots,k\}.$ If we let $k=N$, the corresponding solution will be an optimal solution to the general version of \textbf{[P4]}.
\end{proposition}

For any $k$, using the $\nu_{k-j}^{k}$ values, we can compute the percentage of individuals assigned to the first route in the AABRP relaxation. Let us denote this percentage by $\eta^{k}_{1}=\frac{n_1^*}{N}$.

\begin{lemma}
\label{lem:mon}
$\eta_1^k$ is monotonically increasing for $k\ge3$.
\end{lemma}
The proof of the lemma is presented in Appendix \ref{Ap:Lemma}.

\begin{proposition}
\label{Prop:Limit}
$\eta^{k}_{1}$ converges to a limiting value as $k \rightarrow \infty$ under a linear form of $S(t)$. 
\end{proposition}

The proof is in Appendix \ref{Ap:Limit}. Our computations reveal that the corresponding limiting value of $\eta_1^k$ is approximately 0.867 as shown in Table \ref{tab:r1_ratio}.  Observe that, as indicated in Lemma \ref{lem:mon}, $\eta_1^k$ initially decreases when $k$ is increased from 2 to 3, but is monotonically increasing for $k\ge3$. 

Moreover, as $k$ increases, the values of $1/\nu_1^k$ and $1/\nu_2^k$, also shown in Table \ref{tab:r1_ratio}, appear to converge approximately to 0.133, which is $1-0.867$. These observations suggest a heuristic for solving the uncapacitated AABRP, where within each route, we serve $86.7\%$ of the remaining individuals (i.e., rounded up to the closest integer) that have not yet been served. We refer to this as the GR heuristic for the uncapacitated AABRP.

\begin{table}[H]
\centering
\caption{$\eta_1^k$, $1/\nu_1^k$, and $1/\nu_2^k$ values for a linear $S(t)$, as $k$ changes.}
\begin{tabular}{@{}cccc@{}}
\toprule
$k$ & $\eta^k_1$  & $1/\nu_1^k$ & $1/\nu_2^k$ \\ \hline
2 & 0.872678  &  0.145898 & - \\
3 & 0.866352  &  0.134624 & 0.145898 \\
4 & 0.866758  &  0.133180 & 0.134624 \\
5 & 0.866977  &  0.132989 & 0.133180 \\
6 & 0.867029  &  0.132964 & 0.132989 \\
7 & 0.867038  &  0.132960 & 0.132964 \\
8 & 0.867040  &  0.132960 & 0.132960 \\
9 & 0.867040  &  0.132960 & 0.132960 \\
10 & 0.867040 &  0.132960 & 0.132960 
\\ \bottomrule
\end{tabular} \label{tab:r1_ratio}
\end{table}

\normalsize

\subsection{Capacitated AABRP}

In this subsection, we assume that the vehicle has a capacity of $C$ individuals. The relaxation of the capacitated AABRP is formulated as follows, which differs from \textbf{[P4]} by capacity constraints only.
\begin{eqnarray}
 \textbf{[P5]} \qquad\mbox{Min} & \sum_{i=1}^k n_i \left(\sum_{j=1}^i \sqrt{n_j}\right)   &  \nonumber\\
\mbox{s.t.} & \sum_{i=1}^k n_i = N, & \nonumber \\
& -n_i \le 0,  & i = 1, \ldots, k. \nonumber \\
& n_i \leq C, & i = 1, \ldots, k. \nonumber
\end{eqnarray}
\normalsize
To write the necessary KKT conditions for the problem in \textbf{[P5]}, we use the same KKT multipliers as in \textbf{[P4]}. In addition, we let $\omega_i$ denote a multiplier for the $i^{\text{th}}$ capacity constraint. Then, the solution for the relaxation of the capacitated AABRP must satisfy the following conditions:
\begin{eqnarray}
 & \frac{3\sqrt{n_i}}{2} + \sum_{j=1}^{i-1}\sqrt{n_j} + \frac{\sum_{j=i+1}^{k}n_j}{2\sqrt{n_i}} - \alpha - \lambda_i + \omega_i  = 0,  &  i = 1, \ldots, k, \label{eq:KKTCap}  \\ 
 & \sum_{i=1}^k n_i = N, & \\ 
 & \lambda_i n_i = 0, & i = 1, \ldots, k,\\ 
 &  \omega_i  (n_i - C) = 0, & i = 1, \ldots, k, \\
 & 0 \leq n_i \leq C, & i = 1, \ldots, k, \\ 
 & \lambda_i, \omega_i  \ge 0, & i = 1, \ldots, k. 
\end{eqnarray}

If the optimal solution to the uncapacitated version obeys the capacity constraints, it is also optimal for the capacitated one. This occurs when $\eta^k_1 N\leq C$. Otherwise, we need to incorporate the KKT conditions of the capacitated version. 

Now, suppose that $\eta^k_1 N > C$, implying that the solution for the uncapacitated problem does not satisfy the KKT conditions for the capacitated problem. For this case, assume that in the optimal solution for the capacitated problem, $n_i<C$ for all $i$. Complementary slackness implies that $\omega_i=0$ for all $i$. Therefore, we can remove $\omega_i$ from Equation \eqref{eq:KKTCap} and inequalities $\omega_i\geq 0$ are satisfied. Then, the remaining equations and inequalities become the same as those in the uncapacitated problem. Hence, the KKT solution for the uncapacitated problem should also satisfy these conditions, which leads to a contradiction. Therefore, we cannot have $n_i<C$ for all $i$ in the optimal solution, implying that there is at least one $n_i$ that is equal to $C$. By Proposition \ref{Prop:ordering}, we also know that $n_i \geq n_{i+l}$ for $l>0$. Therefore, if exactly one route is at capacity, it must be the first route, i.e., $n_1=C$. By removing these $C$ individuals from the entire set, we can solve the relaxation of the capacitated AABRP problem for the remaining $N-C$ nodes. Following the same argument, if the solution to the uncapacitated problem for $N-C$ nodes is also feasible for the capacitated problem, it is optimal for the capacitated one. Otherwise, the next route, which is the second route for the original problem, must have $C$ individuals. We repeat the same procedure until we assign all individuals.

The following proposition summarizes an optimal solution for the relaxation of the capacitated AABRP.
\begin{proposition}
An optimal solution for \textbf{[P5]} is shown in Table \ref{Tab:CapOptimal}. If we let $k=N$, the corresponding solution will be optimal for the general version of \textbf{[P5]}.

\renewcommand{\arraystretch}{2}
\begin{table}[h]
\centering
\caption{Optimal solution for the relaxation of the capacitated AABRP.}
\scalefont{0.65}
\begin{tabular}{@{} c| c c c c c  c@{}} 
\toprule
\multirow{2}{*}{Condition} & \multicolumn{6}{c}{Optimal Solution} \\ \cline{2-7} 
& $n_1^*$ & $n_2^*$  & $n_3^*$& $\cdots$ &$n_{k-1}^*$& $n_k^*$ \\ \hline
$1 < \frac{N}{C} \leq \frac{1}{\eta^k_1}$ & $\eta^k_1 N $ & $\frac{\eta^k_1}{\nu^k_1} N $ & $\frac{\eta^k_1}{\nu^k_1 \nu^k_2}  N$ &  $\cdots$ & $\frac{\eta^k_1 }{\prod_{j=1}^{k-2} \nu^k_j} N$ & $\frac{\eta^k_1 }{\prod_{j=1}^{k-1} \nu^k_j} N$ \\
$\frac{1}{\eta^k_1} < \frac{N}{C} \leq 1+\frac{1}{\eta^{k-1}_1}$ & $C$ &  $\eta^{k-1}_1 (N-C) $ &  $\frac{\eta^{k-1}_1}{\nu^{k-1}_1} (N-C) $ & $\cdots$ & $\frac{\eta^{k-1}_1 }{\prod_{j=1}^{k-3} \nu^{k-1}_j} (N-C)$ & $\frac{\eta^{k-1}_1}{\prod_{j=1}^{k-2} \nu^{k-1}_j}  (N-C)$ \\
$1+\frac{1}{\eta^{k-1}_1} < \frac{N}{C} \leq  2+ \frac{1}{\eta^{k-2}_1}$ & $C$ & $C$ &  $\eta^{k-2}_1 (N-2C) $ &  $\cdots$ & $\frac{\eta^{k-2}_1}{\prod_{j=1}^{k-4} \nu^{k-2}_j}  (N-2C)$ & $\frac{\eta^{k-2}_1}{\prod_{j=1}^{k-3} \nu^{k-2}_j}  (N-2C)$ \\
$\vdots$ & $\vdots$ & $\vdots$ & $\vdots$ & $\ddots$ & $\vdots$  & $\vdots$ \\
$k-3+\frac{1}{\eta^3_1}  < \frac{N}{C} \leq k-2+\frac{1}{\eta^2_1}$ & $C$ & $C$ & $C$ & $\cdots$& $\eta^2_1(N-(k-2)C)$ &  $\frac{\eta^2_1}{\nu^2_1}  (N-(k-2)C)$ \\
$k-2+\frac{1}{\eta^2_1}  < \frac{N}{C} < k$ & $C$ & $C$ & $C$ & $\cdots$& $C$ &  $N - (k-1)C$ \\\bottomrule
\end{tabular}
\label{Tab:CapOptimal}
\end{table}

\normalsize
\end{proposition}

The GR heuristic, proposed for the uncapacitated version, can thus be adapted to solve the capacitated AABRP, resulting in the GR heuristic for the capacitated AABRP, presented in Algorithm \ref{Alg:GoldRatio}.

\begin{algorithm}[h]
\caption{GR Heuristic for the Capacitated AABRP}
\label{Alg:GoldRatio}
\begin{algorithmic}[1]
\Require The number of nodes $N$ and the vehicle capacity $C$.
\Ensure The number of individuals on each route, i.e., $n_i$'s, in the AABRP.
\State $\eta \gets  0.867$
\State $\text{Remaining} \gets N$
\State $i\gets 0$
\While{$\text{Remaining} >0$}
\State $i\gets i+1$
\State $n_i = \lceil \min\{C, \eta\cdot \text{Remaining}\}\rceil$
\State $\text{Remaining} \gets \text{Remaining} - n_i$
\EndWhile
\State Return $n_1, n_2, \dots, n_i$
\end{algorithmic}
\end{algorithm}

\section{Computational Experiments}
\label{Sec:Comp}

To evaluate the computational performance of our GR heuristic in solving the AABRP  with a linear satisfaction function $S(t)=a-bt$, in this section we perform numerical experiments and compare our heuristic solutions with the exact solutions. We solve the AABRP problem in formulation [\textbf{P2}] using the BARON solver and the GR heuristic, respectively, and provide the associated cumulative satisfaction value as in objective function \eqref{eq:P3_maxObj}. All tested problems are solved on a Windows 11 desktop with an AMD Ryzen 7 5800X 8-Core processor and 16GB memory. 

To test a sufficient number of instances, we choose a vehicle capacity of $C$ from $\{16, 18, 20\}$, and for each value of $C$ we increase the number of nodes $N$ from 10 to 100 with a step size of 2. In addition to solving the capacitated instances, we also solve a set of uncapacitated instances by using both methods. We vary the value of $N$ in the same manner as for the capacitated instances. With $a=1$, $b=0.01$, and $\kappa=1$, a total of 184 base instances are solved. 
Among these, in 73 instances, our GR heuristic finds the exact same solution as BARON, and in the remaining 111 instances, the average deviation in the objective function is less than $0.02\%$. Moreover, the running times are significantly improved when using the GR heuristic, especially for large values of $N$.

Table \ref{tab:baron_vs_gr} compares the BARON and GR solutions for selected combinations of $C$ and $N$, including the objective values (i.e., the cumulative satisfaction $z$), the number of individuals in each route ($n_i$'s), and the running times of the two methods. The objective values are rounded to four decimal places. Instances are starred if the GR solution and the BARON solution differ. All instances can be solved to optimality in BARON within or around a minute. However, our GR heuristic can obtain a high-quality solution in a negligible amount of computational time. It should be mentioned that we utilize constraints $n_1 \geq n_2 \geq \ldots \geq n_N$, which are validated by Proposition \ref{Prop:ordering}, to expedite BARON. Otherwise, any instance with $N \geq 20$ takes at least one hour to terminate at BARON's optimality, and as $N$ increases, it can take more than five hours. 

\renewcommand{\arraystretch}{1}
\begin{table}[h]
\centering
\caption{Numerical solutions for the AABRP using BARON and the GR heuristic.}
\resizebox{\textwidth}{!}{%
\begin{tabular}{@{}ccl|ccc|ccc@{}}
\toprule
\multicolumn{1}{l}{} & \multicolumn{1}{l}{} &  & \multicolumn{3}{c|}{\textbf{BARON Solution}}   & \multicolumn{3}{c}{\textbf{GR Heuristic Solution}} \\ \hline
$C$ & $N$ &  & $z$ & $[n_i]$ & Time (s)  & $z$  & $[n_i]$ & \multicolumn{1}{l}{Time (s)} \\ \hline
\multirow{3}{*}{16} & 20 & * & 19.1207  & {[}16, 3, 1{]} & 0.69  & 19.1200 & {[}16, 4{]} & \multirow{3}{*}{\textless{}0.01} \\
 & 40 &  & 37.2183 & {[}16, 16, 7, 1{]} & 1.22   & 37.2183 & {[}16, 16, 7, 1{]} &  \\
 & 100 & * & 85.5207 & {[}16, 16, 16, 16, 16, 16, 3, 1{]} & 20.27  & 85.5200 & {[}16, 16, 16, 16, 16, 16, 4{]} &  \\
 \hline
\multirow{3}{*}{18} & 20 & * & 19.1234 & {[}17, 3{]} & 0.40  & 19.1232 & {[}18, 2{]} & \multirow{3}{*}{\textless{}0.01} \\
 & 40 & * & 37.2903 & {[}18, 18, 3, 1\} & 3.32   & 37.2896 & {[}18, 18, 4{]} &  \\
 & 100 &  & 86.1135 & {[}18, 18, 18, 18, 18, 9, 1{]} & 15.80   & 86.1135 & {[}18, 18, 18, 18, 18, 9, 1{]} &  \\
 \hline
\multirow{3}{*}{20} & 20 & * & 19.1234 & {[}17, 3{]} & 0.42 & 19.1232 & {[}18, 2{]} & \multirow{3}{*}{\textless{}0.01} \\
 & 40 & * & 37.3346 & {[}20, 17, 3{]} & 1.59  & 37.3343 & {[}20, 18, 2{]} &  \\
 & 100 & * & 86.6014 & {[}20, 20, 20, 20, 17, 3{]} & 62.27   & 86.6012 & {[}20, 20, 20, 20, 18, 2{]} &  \\
 \hline
 \multirow{3}{*}{Uncap.} & 20 & * & 19.1234 & {[}17, 3{]} & 0.69 & 19.1232 & {[}18, 2{]} & \multirow{3}{*}{\textless{}0.01} \\
 & 40 & * & 37.5236 & {[}35, 4, 1{]} & 2.24  & 37.5218 & {[}35, 5{]} &  \\
 & 100 & * & 90.2132 & {[}87, 11, 2{]} & 42.04   & 90.2123 & {[}87, 12, 1{]} &  \\
 \bottomrule
\end{tabular} \label{tab:baron_vs_gr}
}\end{table}

 When $N=20$ and $C=20$, although one route can accommodate all 20 individuals, both methods provide a solution that requires an additional route, which leads to a higher total satisfaction value than transporting all individuals on a single route. Without the capacity constraints, when $N=100$, all individuals are accommodated using three routes.

\section{Conclusions and Future Work}
\label{Sec:Conc}
In this paper, we study group routing decisions where a number of individuals need to be transported from their original location to a facility, and the time it takes to arrive at the facility has a negative impact on their well-being or decreases their satisfaction. We call the underlying routing problem the ABRP, which is NP-complete. Such problems emerge in humanitarian logistics, disaster response, military operations, and service systems, where finding good solutions in a short amount of time is critically important. Hence, the ABRP poses high-value research challenges.
Motivated by the computational issues, we focus on developing general policies to support dispatchers' decisions in real time using approximation. We develop our approximate model using \cite{beardwood1959shortest}'s seminal TSP tour length estimation model. Assuming that the individuals are uniformly dispersed in the service area and satisfaction is a linear decreasing function of time, we develop some structural results to characterize the optimal solution. We also devise a heuristic based on these structural results. The computational experiments show that our heuristic finds near-optimal solutions with less than a second. 

Our work opens up new research questions. One future research direction is to explore group routing decisions for more general satisfaction functions. In addition, in emergency situations, routing operations are generally followed by other service operations, such as healthcare services in destination facilities. We plan to explore integration of routing decisions with capacity allocation decisions to improve synchronization between these operations and overall system efficiency.

\bibliographystyle{chicago}

\bibliography{AmbulanceRouting}

\appendix
\section{Proof of Proposition \ref{Prop:NP-complete}}
\label{Ap:NPComplete}

\begin{proof}
To prove that the ABRP is NP-complete, we will reduce the shortest Hamiltonian cycle problem, which is known to be NP-complete, to the ABRP. First, we will show how to modify a shortest Hamiltonian cycle problem instance for the ABRP such that the optimal solution of the ABRP has a single route.

Consider an arbitrary instance for the shortest Hamiltonian cycle problem with $n$ nodes. Let $H$ be the length (i.e., duration in our case) of the shortest Hamiltonian path over $n$ nodes. Also, consider $k$ arbitrary subsets of $n$ nodes, with $n_1, n_2, \ldots, n_k$ nodes in each subset and $\sum_{i=1}^{k}n_i=n$.
Let $H_1, H_2, \ldots, H_k$ denote the optimal Hamiltonian path lengths over the subsets $1, 2, \dots, k$.

Suppose that we add a dummy node, which represents the depot, to the current node set, that has a distance $M$ to all existing nodes. Assume that this dummy node is both the depot where the ambulance bus is located and the destination facility for individuals, and assume that the ambulance bus has no capacity limit, i.e., $C=\infty$. We solve the ABRP for this modified instance. In the general case, the optimal solution to the ABRP may have $k \geq 1$ routes. For the optimal ABRP solution to have a single trip, the following conditions must hold for all $k\geq 2$:
\begin{equation}
n(2M+H) < n_1(2M+H_1)+ n_2(4M+H_1+H_2)+\cdots + n_k(2kM+H_1+H_2+\cdots+H_k). \label{Eq:SingleRoute}
\end{equation}
Notice that the right-hand side of this inequality is equal to $2M\sum_{i=1}^{k} n_i + 2n_2M$ plus some other positive terms. Therefore, inequality \eqref{Eq:SingleRoute} holds if the following condition holds:
\begin{equation}
n(2M+H) < 2M\sum_{i=1}^{k} n_i + 2n_2M=2Mn +2n_2M.\label{Eq:SingleRoute1}
\end{equation}
Since $n_2\geq 1$ for $k\geq 2$, by choosing $M=\frac{n}{2}H+\epsilon$, where $\epsilon$ is an infinitesimally small positive quantity, we arrive at a special case of the ABRP where strict inequality in \eqref{Eq:SingleRoute1} holds and the optimal solution to the ABRP has a single trip.

When $k=1$ in the optimal solution for the ABRP, the objective function value is $n$ times the duration of the single route starting from the depot, visiting all nodes, and returning back to the depot. The former (i.e., $n$) is fixed, and the latter can be minimized if the vehicle follows the shortest TSP tour. Since the dummy node has the same distance to all other nodes, by removing the dummy node from the shortest TSP tour, we obtain the shortest Hamiltonian path over the original set of $n$ nodes, which completes the reduction of the shortest Hamiltonian path problem to the ABRP.   \end{proof}

\section{Proof of Proposition \ref{prop:routeorder}}
\label{Ap:routeorder}

\begin{proof}
Consider a solution with route durations $t_1, \ldots, t_{k-1}, t_k, \ldots, t_{h-1}, t_h, \ldots, t_K$, and assume that $n_k \leq n_h$ and $t_k > t_h$. That is, we have a longer duration route containing fewer individuals  that is visited before a shorter duration route with more individuals. Let $z$ be the corresponding objective function value, that is, $z=\sum_{q=1}^{K}n_qS(T_q)$,
where $T_q=\sum_{p=1}^{q}t_p$.

Now consider a different solution where we swap the $k^{\text{th}}$ and $h^{\text{th}}$ routes. By doing so, the duration of each route remains the same, but the arrival times at the destination change for some routes. The new objective function value $z^{\prime}$ is calculated as follows: $$z^{\prime}=\sum_{q=1}^{k-1}n_q S(T_q) + \sum_{q=k}^{h}n_q S(T_q^{\prime}) + \sum_{q=h+1}^{K}n_qS(T_q),$$ where $T_q^{\prime}$ is the new route completion time for $k \leq q \leq h$.

By subtracting $z^{\prime}$ from $z$, we have
\begin{equation}
\begin{split}
z-z^{\prime} = & n_k S(T_{k-1}+t_k) + n_{k+1}S(T_{k-1}+t_k+t_{k+1}) + \cdots \\
& + n_{h-1} S(T_{k-1}+ t_k + \cdots +t_{h-1}) + n_h S(T_{k-1}+ t_k + \cdots +t_h) \\&- \Big(n_hS(T_{k-1}+t_h) + n_{k+1} S(T_{k-1} + t_h + t_{k+1}) + \cdots \\
& + n_{h-1} S(T_{k-1} +t_h+ t_{k+1}  + \cdots +t_{h-1}) + n_k S(T_{k-1} +t_h+ t_{k+1}  + \cdots +t_{h-1} +t_k)\Big).\nonumber
\end{split}
\end{equation}
Since $S(t)$ is a non-increasing function and $t_h < t_k$, the differences between the middle terms are negative. By replacing these terms with $\theta<0$, we can rewrite $z-z^{\prime}$ as follows:
\begin{equation}
\label{Eq:Diff}
\begin{split}
z-z^{\prime}=& n_k \Big{(}S(T_{k-1}+t_k)- S(T_{k-1} +t_h+ t_{k+1} + \cdots +t_{h-1}+t_k)\Big{)} \\
&+ n_h \big{(} S(T_{k-1}+t_k+ \cdots + t_h) -S(T_{k-1}+t_h) \big{)}
+ \theta. 
\end{split}
\end{equation}
Since $n_k\leq n_h$ and $S(t)$ is non-increasing, we can obtain the following inequality by replacing $n_k$ with $n_h$ in Equation (\ref{Eq:Diff}).
\begin{equation}
\begin{split}
z-z^{\prime} \leq & n_h \Big{(}S(T_{k-1}+t_k)- S(T_{k-1} +t_h+ t_{k+1}  + \cdots +t_{h-1}+t_k)\Big{)}\\
& + n_h \Big{(} S(T_{k-1}+t_k+ \cdots + t_h) -S(T_{k-1}+t_h) \Big{)} + \theta \\
= & n_h \Big{(} S(T_{k-1}+t_k) - S(T_{k-1} + t_h) \Big{)} + \theta < 0.  \nonumber
\end{split}
\end{equation}
Therefore, $z< z^{\prime}$, which implies that the initial solution cannot be optimal. 
\end{proof}

\section{Proof of Proposition \ref{Prop:ordering}}
\label{Ap:ordering}
\begin{proof}
For the AABRP with a linear and decreasing satisfaction function, suppose we have optimal values for $n_3$, $n_4$, $\dots, n_k$, such that $\sum_{i=3}^k n_i=M$, $n_1>0$, and $n_2>0$. The objective function terms that contain $n_1$ and $n_2$, denoted by $h(n_1, n_2)$, can be written as follows:
\begin{equation}
h(n_1, n_2) = n_1 \sqrt{n_1} + n_2(\sqrt{n_1} + \sqrt{n_2}) + M (\sqrt{n_1} + \sqrt{n_2}). \nonumber
\end{equation}

Assume that in the optimal solution we have $n_2>n_1$. Consider another solution with $n_2^{\prime}=n_1$ and $n_1^{\prime}=n_2$. For these two solutions, we have
\begin{equation}
    h(n_1^{\prime}, n_2^{\prime}) - h(n_1, n_2) = n_1\sqrt{n_2} - n_2 \sqrt{n_1} = \sqrt{n_1 n_2} (\sqrt{n_1} - \sqrt{n_2}) .\nonumber
\end{equation}

The difference is negative if and only if $\sqrt{n_1} < \sqrt{n_2}$, which is true by assumption. Thus, the solution with $n_1^{\prime}$ and $n_2^{\prime}$ is an improvement over the optimal solution, which is a contradiction. This implies that if an optimal solution satisfies $n_1>0$ and $n_2 >0$, then $n_1\geq n_2$. The same argument can be applied recursively to any pair of variables $n_i$ and $n_{i+1}$. Therefore, in the optimal solution, $n_i$'s are successively non-increasing.
\end{proof}

\section{Proof of Proposition \ref{Prop:Strictpositve}}
\label{Ap:Strictpositve}
\begin{proof}
Suppose that in an optimal solution $ [{n}_i^\star]= ({n}_1^\star,{n}_2^\star,\ldots,{n}_{k-1}^\star,{n}_k^\star)$, we have some $n_i^\star$ which is equal to 0. Due to Proposition \ref{Prop:ordering}, without loss of generality, we assume that ${n}_k^\star = 0$. Now, let us construct another solution $[{n}_i^\circ] = (n_1^\star,n_2^\star,\ldots,n_{k-2}^\star,n_{k-1}^\star-\delta,\delta),$ where $0 < \delta < n_{k-1}^\star.$ It is obvious that $[{n}_i^\circ]$ is feasible. Let $f([{n}_i^\star])$ and $f([{n}_i^\circ])$ be the corresponding objective function values evaluated at $[{n}_i^\star]$ and $[{n}_i^\circ]$, respectively. Therefore,
\begin{eqnarray}
f([{n}_i^\circ]) - f([{n}_i^\star]) & =  & (n_{k-1}^\star-\delta) \left(\sum_{j=1}^{k-2} \sqrt{{n}_j^\star} + \sqrt{{n}_{k-1}^\star - \delta}\right) +  \nonumber \\
& &  \delta \left(\sum_{j=1}^{k-2} \sqrt{{n}_j^\star} + \sqrt{{n}_{k-1}^\star - \delta} + \sqrt{\delta}\right)  -  {n}_{k-1}^\star \sum_{j=1}^{k-1} \sqrt{{n}_j^\star} \nonumber \\
& =  &  {n}_{k-1}^\star \sqrt{{n}_{k-1}^\star - \delta} + \delta \sqrt{\delta} - {n}_{k-1}^\star \sqrt{{n}_{k-1}^\star} \nonumber  
\end{eqnarray}
If we replace ${n}_{k-1}^\star$ by $A$ and let $g(\delta) = A \sqrt{A - \delta} + \delta \sqrt{\delta} - A\sqrt{A}$, we have $g^{\prime}(\delta) = \frac{-A}{2\sqrt{A-\delta}}+\frac{3}{2} \sqrt{\delta}$. Therefore, $g^{\prime}(\delta)|_{\delta = 0} = - \frac{1}{2} \sqrt{A} < 0.$ Moreover, we observe that $g(0) = 0$. Therefore, there exists a $\delta$ within the range of $0$ and $n_{k-1}^\star$ such that $f([{n}_i^\circ]) - f([{n}_i^\star]) < 0,$ which contradicts the optimality of $[{n}_i^\star]$. 
\end{proof}

\section{Proof of Lemma \ref{lem:mon}}
\label{Ap:Lemma}

\begin{proof}
By recursively applying the identity $n_{k-l}= \nu_{k-l}^k n_{k-l+1}$, we can write
\begin{eqnarray}
    \eta_1^k =\frac{\prod_{j=2}^k \nu_1^j}{1+\sum_{i=2}^k \prod_{j=2}^i\nu_1^j}.
\end{eqnarray}
Using the above equation, we can show that $\eta_1^{k+1}-\eta_1^k$ equals a positive multiple of 
\begin{eqnarray}
    \nu_1^{k+1}+\sum_{i=2}^{k-1}(\nu_1^{k+1}-\nu_1^{i+1})\prod_{j=2}^i \nu_1^j - (\nu_1^2+1). \label{eq:posmult}
\end{eqnarray}
Using Equation \eqref{Eq:NUs}, we can show that $\nu_1^j$ is strictly increasing in $j$ with $\nu_{1}^2=(1+\varphi)^2=\left( \frac{3+\sqrt{5}}{2}\right)^2\approx 6.854$.  This implies that $\nu_1^{k+1}$ is positive and increasing in $k$, as is the second term in Expression \eqref{eq:posmult} (i.e., $\sum_{i=2}^{k-1}(\nu_1^{k+1}-\nu_1^{i+1})\prod_{j=2}^i \nu_1^j$) for $k\ge3$, while the third term, $\nu_1^2+1$, is fixed for any $k$ (at approximately 7.854). When $k=3$, the first term in Equation \eqref{eq:posmult} equals approximately 7.509 and the second term equals approximately 0.552, for a total of approximately 8.061.  Thus $\eta_1^4-\eta_1^3$ is a positive multiple of $0.207$ and is greater than 0.  Because the first two terms in Equation \eqref{eq:posmult} are positive and increasing in $k$ and the third term is fixed, Equation \eqref{eq:posmult} is positive for any $k\ge 3$, which implies that $\eta_1^k$ is monotonically increasing in $k$. It is worth noting that when $k=2$, Expression \eqref{eq:posmult} becomes $\nu_1^{3}- (\nu_1^2+1)\approx 7.428-7.854<0$, and the monotonicity property only holds for $k\ge 3$.  
\end{proof}

\section{Proof of Proposition \ref{Prop:Limit}}
\label{Ap:Limit}
\begin{proof}
From  Equation \eqref{eq:n1}, we have
\begin{eqnarray}
    \eta_1^k = \frac{n_1}{N}=\frac{1}{1+\sum_{i=1}^{k-1}\hat{\rho}_i^k}=\frac{1}{1+\frac{1}{\nu_1^k}+\frac{1}{\nu_1^k\nu_2^k}+\cdots + \frac{1}{\nu_1^k\nu_2^k\cdots \nu_{k-1}^k}}.\label{eq:eta1k}
\end{eqnarray}  
Consider the denominator of the above expression, which is bounded from below by (and strictly greater than) one, implying that $\eta_1^k$ is bounded from above by (and strictly less than) one.  Next, observe that the denominator in \eqref{eq:eta1k} is strictly less than  
\begin{eqnarray}
    1+\frac{1}{\nu_{k-1}^k}+\frac{1}{(\nu_{k-1}^k)^2}+\frac{1}{(\nu_{k-1}^k)^3}+\cdots + \frac{1}{(\nu_{k-1}^k)^{k-1}},\label{eq:denombound}
\end{eqnarray} 
which we can write as
\begin{eqnarray}
    1+r+r^2+r^3+\cdots + r^{k-1},\label{eq:denombound}
\end{eqnarray} 
with $r=\frac{1}{\nu_{k-1}^k}\approx\frac{1}{6.854}\approx0.146$.  Thus, \eqref{eq:denombound} is a convergent sequence as $k\rightarrow \infty$ with limit approximately $\frac{1}{1-r}=\frac{1}{1-\frac{1}{\nu_{k-1}^k}}\approx 1.171$.  Because the denominator of \eqref{eq:eta1k} is bounded above by \eqref{eq:denombound}, this implies that \eqref{eq:eta1k} is bounded from below by $\frac{1}{1.171}\approx 0.854$, i.e., \eqref{eq:eta1k} falls between 0.854 and 1. Because $\eta_1^k$ is monotonically increasing in $k$ by Lemma \ref{lem:mon} and is bounded between 0.854 and 1, this implies that $\eta_1^k$ converges to a value in this interval as $k\rightarrow \infty$.  
\end{proof}
\end{document}